\documentstyle{amsart}

\newtheorem{theorem}{Theorem}[section]
\newtheorem{defn}[theorem]{Definition}
\newtheorem{prop}[theorem]{Proposition}
\newtheorem{lem}[theorem]{Lemma}

\newtheorem{claim}[theorem]{Claim}

\newcommand{\sub}{\subseteq}

\newcommand{\rt}{\rightarrow}
\newcommand{\sm}{\setminus}
\newcommand{\proof}{\noi {\bf Proof:}\ }
\newcommand{\noi}{\noindent}

\newsymbol\upharpoonright 1316
\newcommand{\rest}{\upharpoonright}

\newcommand{\ra}{\rangle}
\newcommand{\la}{\langle}

\newcommand{\frt}{{\frak t}}
\newcommand{\frf}{{\frak f}}

\newcommand{\bR}{{\Bbb R}}

\newcommand{\CA}{{\cal A}}

\newcommand{\CC}{{\cal C}}
\newcommand{\CF}{{\cal F}}
\newcommand{\CG}{{\cal G}}
\newcommand{\CH}{{\cal H}}
\newcommand{\CI}{{\cal I}}
\newcommand{\CN}{{\cal N}}
\newcommand{\CO}{{\cal O}}

\newcommand{\CW}{{\cal W}}   
\newcommand{\CX}{{\cal X}}

\newcommand{\Fr}{{\frak Fr }}
\newcommand{\Fs}{F$_\sigma\;\;$}
\newcommand{\dom}{dom}

\newcommand{\fo}{\; ^\omega  \! \omega}

\newcommand{\fso}{[\omega]^{<\omega}}
\newcommand{\ff}{\:^{\omega} \! 2}
\newcommand{\fff}{\:^{<\omega} \! 2}

\newcommand{\iso}{[\omega]^\omega}
\newcommand{\po}{\mbox{{\Large $\wp$}}(\omega)}
\newcommand{\fa}{\forall}					  
\newcommand{\fai}{\forall^\infty}					  
\newcommand{\ex}{\exists}
\newcommand{\exi}{\exists^\infty}
\newcommand{\lbv}{[\![}
\newcommand{\rbv}{]\!]}

\begin{document}

\title{Combinatorial aspects of F$_\sigma$ filters \\
   with an application to $\CN$-sets}

\author{Claude Laflamme}
\address{Department of Mathematics and Statistics \\
         University of Calgary \\
         Calgary, Alberta  \\
         Canada T2N 1N4}
\email{laflamme@@acs.ucalgary.ca}


\date{}

\thanks{This research was partially supported by  NSERC of Canada.}

 \maketitle

\begin{abstract}

We discuss F$_\sigma$ filters and show that the minimum size of a
filter base generating an undiagonalizable filter included in some
F$_\sigma$ filter is the better known bounded evasion number ${\frak
e}_{ubd}$. An application to $\CN$-sets from trigonometric series is
given by showing that if $A$ is an $\CN$-set and $B$ has size less
${\frak e}_{ubd}$, then $A \cup B$ is again an $\CN$-set.

\end{abstract}
							  
\section{Introduction}

Our terminology is standard but we review the main concepts and
notation.  The natural numbers will be denoted by $\omega$, $\po$
denotes the collection of all its subsets.  Given a set $X$, we write
$[X]^{\omega}$ and $[X]^{<\omega}$ to denote the infinite or finite
subsets of $X$ respectively; if we wish to be more specific, we write
$[X]^{n}$ and $[X]^{\leq n}$ to denote subsets of size $n$ or at most
$n$ respectively.  We use the well known `almost inclusion' ordering
between members of $\iso$, i.e.  $X \sub^*Y$ if $X \sm Y$ is finite.
We identify $\po$ with $\ff$ via characteristic functions.  The space
$\ff$ is further equipped with the product topology of the discrete
space $\{0,1\}$; a basic neighbourhood is then given by sets of the
form

\[ \CO_s = \{f \in \ff: s \sub f \} \] 

\noi where $s \in \fff$, the collection of finite binary sequences.  The
terms ``nowhere dense'', ``meager'', ``Baire property'' all refer to
this topology. 

A filter is a collection of subsets of $\omega$ closed under finite
intersections, supersets and containing all cofinite sets; it is
called proper if it contains only infinite sets; thus the collection
of cofinite sets is the smallest proper filter, it is called the {\em
Fr\'echet} filter and is denoted by $\Fr$. To avoid trivialities, we
shall assume that all filters under discussion are proper. An infinite
set $X \in \iso$ is said to diagonalize a filter $\CF$ if $X \sub^* Y$
for each $Y \in \CF$.  Given a collection of sets $\CX \sub \iso$, we
denote by $\la \CX \ra$ the filter generated by $\CX$, that is the
smallest filter containing each member of $\CX$.

The Rudin-Keisler ordering on filters is defined by

\[ \CF \leq_{RK} \CG \mbox{ if } (\ex f \in \fo) \; \CG \supseteq
 \{f^{-1}\{X\} : X \in \CF \}.  \]

The following Lemma from \cite{L1} combinatorially describes \Fs
filters.

\begin{lem}
Let $\CF$ be an \Fs filter and $g \in \fo$. Then there is an
increasing sequence of natural numbers $\la n_k : k \in \omega \ra$
and sets $a^k_i \sub [n_k, n_{k+1})$, $i<m_k$, such that
\begin{enumerate}
\item $(\fa x \in [m_k]^{\leq g(k)}) \; \bigcap_{i \in x} a^k_i \neq
\emptyset$ \\
\item $(\fa X \in \CF)(\fai k)( \ex i<m_k) \; a^k_i \sub X. $
\end{enumerate}
\end{lem}

\proof Let $\CF = \cup_n \CC_n$ where each $\CC_n$ is closed and put
$\CC = \{X \cup n: n \in \omega \mbox{ and } X \in \CC_n \}$. Then
again $\CC$ is a closed set and every member of $\CF$ is almost equal
to a member of $\CC$. 

\noi Let $n_0 = 0$ and having defined $n_j$ together with the $a^j_i$'s
for $j \leq k$, choose an $n_{k+1} > n_k$ such that
\[(\fa X_0,X_1,\cdots,X_{g(k)-1} \in \CC) \; \bigcap_{i<g(k)} 
  X_i \cap [n_k,n_{k+1}) \neq \emptyset. \]

\noi The existence of such an $n_{k+1}$ follows from the fact that
$\CC$ is closed and that $\CF$ only contains infinite sets. Now
enumerate $\{X \cap [n_k, n_{k+1}): X \in \CC \}$ as $\{a^k_i : i <
m_k \}$ and this completes the proof. \qed

\bigskip

It is worth noticing that conversely, given a family $\la \la a^k_i :
i <m_k, \ra; k \in \omega , g \ra$ satisfying conditions (1) and (2)
above, then the  collection

\[ \{X: (\fa k)(\ex i < m_k) \; a^k_i \sub X \} \]

\noi is a closed set  generating an \Fs filter whenever $\lim_n g(n) =
\infty$.

\section{Non trivial F$_\sigma \:$ filters}

We generally define a non trivial filter as one that cannot be
diagonalized by a single infinite set. We first present a
combinatorial description of the smallest size of family of sets
generating a non trivial filter included in some F$_\sigma$ filter.
This is a variation of some well known cardinals; indeed the cardinal
$\frak p$ is defined as the smallest size of a family of sets
generating a non trivial filter and $\frak t$ is defined as the
smallest size of well ordered (under almost inclusion) family of sets
generating a non trivial filter. It turns out that these cardinals
have a substantial impact on the Set theory of the Reals.

\begin{defn}\hfill

\begin{tabbing}
${\frak f}  =$ \= $\min \{|\CX| :$ 
 $\CX$ generates a  non trivial filter 
                 included in some F$_\sigma \:$ filter $\}$.\\

   \\

${\frak f}_1 =  \min \{|\CH|:$
     $ \CH \sub \fo$ is bounded and for some $ g \in \fo$
        with $\lim_{n \rightarrow \infty}g(n)=\infty$, \\
    \>$(\fa X \in \iso)(\fa s_n \in [\omega]^{\leq g(n)})\; 
        (\ex h \in \CH)(\exi n \in X) \; h(n) \notin s_n.\} $ \\
\\

${\frak f}_2 =  \min \{|\CH|:$
  $ \CH \sub \fo$ is bounded and for some $ g \in \fo$
        with $\lim_{n \rightarrow \infty}g(n)=\infty$, \\
 \> $(\fa X \in \iso)(\fa \pi_n :\:^n\!\omega \rightarrow
[\omega]^{\leq g(n)})\; 
 (\ex h \in \CH)(\exi n \in X) \; h(n) \notin \pi_n(h \rest n).\} $ \\

 \\

${\frak e}_{ubd} =  \min \{|\CH|:$  $ \CH \sub \fo$ is bounded and  \\
 \>$(\fa X \in \iso)(\fa \pi_n :\: ^n\!\omega \rightarrow \omega)\; 
   (\ex h \in \CH)(\exi n \in X) \; h(n) \neq \pi_n(h \rest n).\} $
\end{tabbing}

\end{defn}

T. Eisworth has shown (unpublished) that ${\frak e}_{ubd} \leq {\frak
f}$ and an argument very similar to that of A. Blass in \cite{Bl}
shows that ${\frak f}_1 \leq {\frak e}_{ubd}$. We extend these results
by showing that all four cardinals are equal.

\medskip

\begin{prop}
The four cardinals ${\frak f}$, ${\frak f}_1$, ${\frak f}_2$ and
${\frak e}_{ubd}$ are equal. 
\end{prop}

\proof For simplicity, we prove $\;{\frak f} \leq {\frak f}_1 \leq {\frak
e}_{ubd} \leq {\frak f}_2 \leq {\frak f}_1 \leq {\frak f}$.

\medskip

\noi {\boldmath ${\frak f} \leq {\frak f}_1$:} Let $\CH \sub \fo$ be
given of size $|\CH| < {\frak f}$, without loss of generality bounded
everywhere by $b \in \fo$ and fix $g \in \fo$ such that $\lim_n g(n) =
\infty$.

\noi Consider $\CA_n = [b(n)]^{\leq g(n)}$ and for $i<b(n)$, put
$a^n_i = \{x \in  \CA_n: i \in x \}$. Notice that 

\[ (\fa x \in  \CA_n=[b(n)]^{\leq g(n)}) 
            \;\; \bigcap_{i \in x} a^n_i \neq \emptyset. \]

\noi Identify $\bigcup_n \CA_n$ with $\omega$, and form the filter
$\CF$ generated by

\[ \{ \bigcup_n a^n_{h(n)}: h \in \CH \}. \]

\noi Then $\CF$ is a filter generated by fewer than $\frak f$ sets and
included in the F$_\sigma$ filter $\la \la a^n_i:i < b(n)\ra: n \in
\omega,g \ra$. Therefore $\CF$ must be trivial, that is for some $X \in
\iso$ and $x_n \in \CA_n$ for $n \in X$, $\{x_n : n \in X \}$
diagonalizes $\CF$. In particular, for $h \in \CH$, 

\[ \{x_n : n \in X \} \sub^* \{a^n_{h(n)} : n \in \omega \}, \]

\noi and thus $h(n) \in x_n$ for all but finitely many $n \in X$.

\medskip

\noi {\boldmath ${\frak f}_1 \leq {\frak e}_{ubd}$:} Let $\CH \sub
\fo$ be given of size $|\CH| < {\frak f}_1$, without loss of
generality bounded everywhere by $b \in \fo$. Partition $\omega$ into
consecutive intervals $\la \CI_n= [a_n, a_{n+1}): n \in \omega \ra$
such that $ a_{n+1} - a_n > n^2$.

\noi For $h \in \CH$, define $\tilde{h}(n)=h \rest \CI_n$ and form
$\tilde{\CH} = \{\tilde{h}: h \in \CH \}$.

\noi Identifying $\prod_{a_n \leq i< a_{n+1}} b(i)$ with its cardinality,
$\tilde{\CH}$ is a bounded family of size  $|\tilde{\CH}| < {\frak
f}_1$. Therefore there is an $\tilde{X} \in \iso$ and $s_n \in 
[ \; \prod_{a_n \leq i< a_{n+1}} b(i) ]^{\leq n}$ such that

\[ (\fa \tilde{h} \in \tilde{\CH})(\fai n \in \tilde{X})\;\;
             \tilde{h}(n) \in s_n. \]

\noi Now by the pigeonhole principle, there must be for each $n \in
\tilde{X}$ an $i_n \in \CI_n$ such that

\[ (*) \hspace{.5in} (\fa t,t' \in s_n)\;\; t \rest i_n = t' \rest i_n \rightarrow
t(i_n) = t'(i_n); \]

\noi this is where we use the fact that $|\CI_n| > n^2$ while $|s_n| \leq
n$. 

\noi Let $X = \{i_n : n \in \tilde{X}\}$ and define
$\pi_i:\:^i\!\omega \rightarrow \omega$ as follows.  If $i=i_n \in X$
and $t \in\: ^i\!\omega$ is such that $t \rest [a_n, i) $ is an initial
segment of a member $t'$ of $s_n$, then define $\pi_i(t)= t'(i_n)$;
this is well defined by the previous paragraph. In all other cases
define $\pi_i(t)$ arbitrarily.

\noi Now for $h \in \CH$, $\tilde{h} \in s_n$ for all but finitely
many $n \in \tilde{X}$; for each such $n$, $i = i_n \in X$, $h \rest
[a_n,i)$ is an initial segment of a member of $s_n$, namely
$\tilde{h}(n) = h \rest \CI_n$, and thus $ \pi_i(h \rest i) = h(i)$.
This proves that ${\frak f}_1 \leq {\frak e}_{ubd}$ as desired.

\medskip

\noi {\boldmath ${\frak e}_{ubd} \leq {\frak f}_2$:} This inequality is
trivial. 

\medskip

\noi {\boldmath${\frak f}_2 \leq {\frak f}_1$:} Let $\CH \sub \fo$ be
given of size $|\CH| < {\frak f}_2$, without loss of generality
bounded everywhere by $b \in \fo$, and fix $g \in \fo$ such that
$\lim_n g(n) = \infty$.

\noi Choose integers $\delta_0=0 < \delta_1 < \ldots$ such that

\[ (\fa n) \; \; \prod_{i \leq \delta_n} b(i) \times g(\delta_n) \leq
g(\delta_{n+1}). \]

\noi Now for $n \in \omega$, define

\[ \tilde{b}(n) = \prod_{\delta_n \leq i \leq \delta_{n+1}} b(i), \]

\noi which we identify with the cartesian product. For $h \in \CH$, define

\[ \tilde{h}(n) = h \rest[\delta_n, \delta_{n+1}] \in \tilde{b}(n) \]

\noi and put $\tilde{\CH} = \{\tilde{h}; h \in \CH \}$, a bounded
family of size less that ${\frak f}_2$. Therefore 

\[ (\ex \tilde{X} \in \iso)(\ex \pi_n:\: ^n\!\omega \rightarrow
[\tilde{b}(n)]^{\leq g(\delta_n)})(\fa \tilde{h} \in \tilde{\CH})(\fai
n \in \tilde{X}) \;\; \tilde{h}(n) \in \pi_n(\tilde{h}\rest n).\]

\noi For $n \in \tilde{X}$, let 

\[ s_{\delta_{n+1}} =    \{ u(\delta_{n+1}): 
  u \in  \pi_n(t) \mbox{ for } t \in \prod_{i<n} \tilde{b}(i) \}. \]

\noi  Then 

\[ |s_{\delta_{n+1}}| \leq \prod_{i \leq \delta_n} b(i) \times
g(\delta_n) \leq g(\delta_{n+1}). \]

\noi and so $s_{\delta_{n+1}} \in [b(\delta_{n+1})]^{\leq
g(\delta_{n+1})}$.
\noi Finally, given $h \in \CH$, and thus $\tilde{h} \in \tilde{\CH}$,

\[ (\fai n \in \tilde{X}) \; \; \tilde{h}(n) \in \pi_n(\tilde{h}\rest
n), \]

\noi and so 

\[ (\fai n \in \tilde{X}) \; \; h(\delta_{n+1}) \in s_{\delta_{n+1}}. \]

\noi We must therefore conclude that ${\frak f}_2 \leq {\frak f}_1$.

\medskip 

\noi{\boldmath ${\frak f}_1 \leq {\frak f}$:}   Let $\CF$ be a
filter generated by $\la A_\alpha : \alpha < \kappa \ra$, $\kappa <
{\frak f}_1$, and included in the F$_\sigma$ filter $\la \la a^k_i:i <
m_k\ra: k \in \omega,g \ra$ where without loss of generality $g(n)=n$.

\noi For each $\alpha < \kappa$, define a function $f_\alpha \in \fo$
such that $a^k_{f_\alpha(k)} \sub A_\alpha$.

\noi Then $\{f_\alpha : \alpha < \kappa \}$ is a bounded family of
size less that ${\frak f}_1$, and therefore 

\[ (\ex X \in \iso)(\ex s_k \in [m_k]^{\leq k})(\fa \alpha)(\fai k \in
X) \; \; f_\alpha(k) \in s_k. \]

\noi We conclude that $ \bigcup_{k \in X} \bigcap_{i \in s_k} a^k_i$
diagonalizes the filter $\CF$ and is therefore trivial. \qed
 
\bigskip

We conclude this section by giving a small perspective on these
cardinals (see \cite{V} for a description of undefined cardinals).  If
one removes the bounded restriction on $\CH$, the cardinal ${\frak
e}_{ubd}$ becomes the evasion number, known as $\frak e$
(\cite{Bl}). Clearly ${\frak e} \leq {\frak e}_{ubd}$ and it has been
proved consistent by Shelah that ${\frak e} < {\frak e}_{ubd}$.
Without the bounded restriction and with $X=\omega$, the cardinal
${\frak f}_1$ is the additivity of measure, $add(\CN)$, (\cite{B}); I
do not know whether removing the bounded restriction alone still
yields the same cardinal.

\noi In any case we have 
\[ add(\CN) \leq {\frak e} \leq {\frak e}_{ubd} = {\frak f},\]

\noi and it is straightforward to show that ${\frak f}$ (through
${\frak f}_1$) is less than or equal to the uniformity of the null and
meager ideals.

\noi On the other hand, we have proved in \cite{L1} the consistency of
${\frak d} <{\frak f}$ (see also \cite{S2}, \cite{Br2}); therefore we must
conclude using the known independence results regarding the Cicho\'n
diagram that ${\frak e}_{ubd}$ is not provably equal to any of the standard
cardinals ${\frak b}, {\frak d}, {\frak t}$ or additivity, uniformity,
cofinality and covering of the null or meager ideals. Of course the
four cardinals $\frak f$, ${\frak f}_1$, ${\frak f}_2$ and ${\frak
e}_{ubd}$ are all equal and we use the more familiar ${\frak e}_{ubd}$
from  now on.

A better provable lower bound for the number ${\frak e}_{ubd}$ (and
thus $\frak f$)is ${\frak t}$; the idea of the proof is from
\cite{BS}. We shall use in the proof the cardinal $\frak b$, the
minimum size of an unbounded family in $\fo$, and the inequality
${\frak t} \leq {\frak b}$ (see \cite{V}).

\medskip

\begin{prop}
${\frak t } \leq {\frak e}_{ubd} \leq 2^{\aleph_0}$.
\end{prop}

\proof We show that ${\frak t} \leq {\frak f}_1$.  Let $\CH = \la
 h_\alpha : \alpha < \kappa \ra$ be a family of size $\kappa < {\frak
 t}$, bounded everywhere by $b \in \fo$, and fix $g \in \fo$ such that
 $\lim_n g(n) = \infty$.

\noi We construct a sequence $\la \phi_\alpha: \alpha \leq \kappa
\ra$ such that:
\begin{enumerate}
\item $\phi_\alpha:\dom(\phi_\alpha)\rightarrow \fso$
 \begin{itemize}
    \item $\dom(\phi_\alpha) \in \iso$
    \item $\phi_\alpha(k) \in [b(k)]^{\leq g(k)}$
    \item $\lim_{k \in \dom(\phi_\alpha)} g(k) - |\phi_\alpha(k)| = + \infty$
 \end{itemize}
\item $(\fa \beta < \alpha)\; \dom(\phi_\alpha) \sub^* \dom(\phi_\beta)$
      and $(\fai k \in \dom(\phi_\alpha))\; \phi_\beta(k) \sub
\phi_\alpha(k)$
\item $(\fa \beta < \alpha)(\fai k \in \dom(\phi_\alpha)) \;
        h_\beta(k) \in \phi_\alpha(k). $
\end{enumerate}

\noi Once we have obtained $\phi_\kappa$, then clearly $s_k =
 \phi_\kappa(k)$ for $ k \in \dom(\phi_\kappa)$ is as desired.

Now to construct the sequence, assume that we already have
$\{\phi_\beta : \beta < \alpha \}$ for some $\alpha \leq \kappa$. 

\noi If $\alpha = \beta + 1$ is a successor ordinal, define
$\phi_\alpha(k)=\phi_\beta(k) \cup \{h_\alpha(k)\}$ for $k \in
\dom(\phi_\beta)$.

\noi For $\alpha$ a limit ordinal, first choose a $\tilde{g}  \in
\fo$ such that  $\lim_k g(k) -\tilde{g}(k) = \infty$ and

\[(\fa \beta < \alpha) (\fai k \in \dom(\phi_\beta)) \;\; |\phi_\beta(k)|
\leq \tilde{g}(k) \]

 \noi This is possible as $\alpha \leq \kappa < {\frak t} \leq {\frak
b}$.

\noi Now for $\beta < \alpha$ define 

\[A_\beta = \{\la k,x \ra : k \in \dom(\phi_\beta) \mbox{ and }
       \phi_\beta(k) \sub x \sub b(k) \mbox{ and } |x| \leq \tilde{g}(k) \}. \]

\noi Clearly $\beta < \gamma < \alpha \implies A_\gamma \sub^*A_\beta$
and, identifying $A_0$ with $\omega$ and as $\alpha < {\frak t}$, we
can find some infinite $A_\alpha \sub^* A_\beta$ for each $\beta <
\alpha$. Finally we put

\[ \phi_\alpha(k) = \mbox{ any $x$ such that $\la k,x \ra \in
A_\alpha$ } \]

\noi and undefined if there is no such $x$. Clearly $\phi_\alpha$ is
as desired. \qed

\bigskip

\section{Trigonometric series and $\CN$-sets}

In this section we prove that if $A$ is an $\CN$-set and $|B|<{\frak
e}_{ubd}$, then $A \cup B$ is also an $\CN$-set. A similar result
was earlier proved with ${\frak p}$ in the role of ${\frak e}_{ubd}$
in \cite{BB} from which we modeled our proof, and then with ${\frak
t}$ in \cite{BS} from which we modeled the above proof of $\frt \leq
{\frak e}_{ubd}$.

\noi It is known that the collection of $\CN$-sets is not in general
closed under unions; in \cite{Ba} two $\CN$-sets of cardinality $\frak
c$ are constructed whose union is not an $\CN$-set. 

\begin{defn}
A set $A \sub \bR$ is called an $\CN$-set (\cite{Ba}) if there is a
sequence of non negative reals $\la a_n:n \in \omega \ra$ such that:
\[ \begin{array}{ll}
(1) & \sum_{n=0}^\infty a_n = + \infty \\
(2) & (\fa a  \in A) \; \sum_{n=0}^\infty a_n | \sin \pi n a | <
\infty .
\end{array} \] 
\end{defn}

\begin{prop}
If $A \sub \bR$ is an $\CN$-set and $|B|<{\frak e}_{ubd}$,
then $A \cup B$ is also an $\CN$-set.
\end{prop}

\proof Fix a sequence of non negative reals $\la a_n:n \in \omega \ra$
as in the definition for the $\CN$-set $A$. As is now standard
procedure (see \cite{BB} and \cite{BS}), we put $s_n = \sum_{i=0}^n
a_i$ and $b_n = a_n / s_n$; then again $\sum_{n=0}^\infty b_n = +
\infty$ and thus we can choose an increasing sequence of natural
numbers $\la \pi_n : n \in \omega \ra$ such that

\[ (\fa k) \; \sum_{i=\pi_k}^{\pi_{k+1}-1} b_i \geq 1. \]

\noi  Now find an unbounded, non decreasing sequence of natural
numbers $\la q_n :n \in \omega \ra$ such that

\[ \sum_{n=0}^\infty a_n/ s_n^{1 + \frac{1}{q_n} } < \infty \]

\noi and let $\epsilon_n = s_n^{-1/q_n}$. 

\noi Finally for $T \sub B$ and $m \in \omega$, define

\[ a^m_T = \{ k\in \omega : 0 \leq k \leq s_m \mbox{ and }
           (\fa t \in T) \; | \sin \pi k m  t | \leq 2 \pi \epsilon_m
\}, \]

\noi and for $|T| \leq q_{\pi_n}/n$, $m\geq \pi_n$, put

\[ b^m_T = \{ a^m_{T'} : T \sub T' \mbox{ and } |T'| \leq q_{\pi_n} \}.  \]

The following claim is from standard number theory and follows from
the pigeon hole principle.

\begin{claim}
If $T \sub B$ and  $|T| \leq q_{\pi_n}$, then $a^m_T \neq \emptyset$
for any $m \geq \pi_n$.
\end{claim}

\proof Fix $n,m$ and let $T= \{t_i : i < |T| \}$. Define
a map
 
\[ c:(1/\epsilon_m)^{|T|}+1 \rt \prod_{i \in |T|} 1/\epsilon_m \]

\noi  by
$c(j) = (\ell_0,\cdots,\ell_{|T|-1})$ if

\[ (\fa i < |T|)(\ex p \in \omega) \; p + \ell_i \epsilon_m 
               \leq j m t_i < p + (\ell_i + 1)\epsilon_m .\]

\noi There must then be two integers $j_1 < j_2
<(1/\epsilon_m)^{|T|}+1$ such that $c(j_1) = c(j_2)$ and let $k = j_2-
j_1$.

\noi Then $k \leq (1/\epsilon_m)^{|T|} \leq (1/\epsilon_m)^{q_{\pi_n}}
\leq (1/\epsilon_m)^{q_m} = s_m$. Now using $\lbv x \rbv$ to denote
the distance from $x$ to the nearest integer, we have, for $t \in T$,

\[ |\sin \pi k m t | \leq \pi \lbv kmt \rbv \leq 2 \pi \epsilon_m \]

\noi as desired. \qed

\noi Therefore $a^m_{T'} \neq \emptyset$ for each $a^m_{T'} \in b^m_T$
and of course, if $T_i \sub B$, $i<n$ and each $|T_i| \leq
q_{\pi_n}/n$, then $ \bigcap_{i<n} b^m_{T_i} \neq \emptyset $ whenever
$m \geq \pi_n$.

Now consider the set
\[ \CW = \{ \la a^m_{T_m} : \pi_n \leq m < \pi_{n+1} \ra :
  n \in \omega \mbox{ and } |T_m| \leq q_{\pi_n} \}, \]

\noi which we may identify with $\omega$, and the filter
 
\[ \CG = \la  \{ X \sub \CW : (\fa n)(\ex T_n \sub B) \; |T_n| \leq
q_{\pi_n}/n \mbox{ and } \prod_{\pi_n \leq m < \pi_{n+1}} b^m_{T_n} \sub
X \} \ra. \]

\noi Then $\CG$ is an F$_\sigma \:$ filter as it is generated by a closed
set and  contains the filter

\begin{eqnarray*}
\lefteqn{ \CF = \la  \{ X \sub \CW : (\ex T \sub B) \; |T| < \infty}  \\
 & & \mbox{ and for all $n$ such that $q_{\pi_n}/n \geq |T|$, }
\prod_{\pi_n \leq m < \pi_{n+1}} b^m_T \sub X \} \ra. 
\end{eqnarray*}

\noi As $\CF$ is generated by fewer than ${\frak e}_{ubd} = \frf$
sets, it must be diagonalized by an infinite set $X\sub \CW$ which, without
loss of generality, is of the form

\[ X = \{ \la a^m_{T_m} : \pi_{n_\ell} \leq m < \pi_{n_{\ell+1}-1} \ra : \ell
\in \omega \} \]

\noi where $|T_m| \leq q_{\pi_{n_\ell}}$ and $n_\ell <
n_{\ell+1}$. For each $\ell$ and $ \pi_{n_\ell} \leq m <
\pi_{n_{\ell+1}-1}$, pick $k_m \in a^m_{T_m}$.

\noi Now clearly $\sum_\ell \sum_{ \pi_{n_\ell} \leq m < \pi_{n_{\ell+1}-1}}
b_m = \infty$ and it remains to show  that

\[ \sum_\ell \sum_{ \pi_{n_\ell} \leq m < \pi_{n_{\ell+1}-1} } b_m |\sin m
k_m \pi x | < \infty \]

\noi for each $x \in A \cup B$.

\noi For $x \in A$,
\[b_m |\sin m k_m \pi x | \leq  b_m k_m | \sin m \pi x|
  \leq  b_m s_m | \sin m \pi x| =  a_m | \sin m \pi x|. \]

\noi Finally for $x \in B$, let $T = \{x\}$ and thus for all but finitely
many $\ell$, for all $ \pi_{n_\ell} \leq m < \pi_{n_{\ell+1}-1} $, 
\[ | \sin m k_m \pi x | \leq 2 \pi \epsilon_m \]
\noi and therefore

\[ b_m | \sin m k_m \pi x | \leq 2 \pi a_m / s_m^{1 + 1/q_m}. \]

\noi This completes the proof. \qed

\bibliographystyle{amsplain}

\end{document}